# Stochastic Optimal Power Flow with Network Reconfiguration: Congestion Management and Facilitating Grid Integration of Renewables

Xingpeng Li
Department of Electrical and Computer Engineering
University of Houston
Email: xli82@uh.edu

Qianxue Xia
School of Electrical and Computer Engineering
Georgia Institute of Technology
Email: qxia31@gatech.edu

*Abstract*— There has been a significant growth of variable renewable generation in the power grid today. However, the industry still uses deterministic optimization to model and solve the optimal power flow (OPF) problem for real-time generation dispatch that ignores the uncertainty associated with intermittent renewable power. Thus, it is necessary to study stochastic OPF (SOPF) that can better handle uncertainty since SOPF is able to consider the probabilistic forecasting information of intermittent renewables. Transmission network congestion is one of the main reasons for renewable energy curtailment. Prior efforts in the literature show that utilizing transmission network reconfiguration can relieve congestion and resolve congestion-induced issues. This paper enhances SOPF by incorporating network reconfiguration into the dispatch model. Numerical simulations show that renewable curtailment can be avoided with the proposed network reconfiguration scheme that relieves transmission congestion in post-contingency situations. It is also shown that network reconfiguration can substantially reduce congestion cost, especially the contingency-case congestion cost.

## NOMENCLATURE

*Sets*:
$C$     Contingencies.
$G$     Online traditional units.
$G(n)$     Online traditional units at bus $n$.
$IR$     Intermittent renewable units.
$IR(n)$     Intermittent renewable units at bus $n$.
$K$     Branches.
$K(n-)$     Branches of which bus $n$ is the from-bus.
$K(n+)$     Branches of which bus $n$ is the to-bus.
$N$     Buses.
$S$     Scenarios.

*Parameters*:
$P_{g,min}$     Minimum output of traditional unit $g$.
$P_{g,max}$     Maximum output of traditional unit $g$.
$P_{g0}$     Initial output of traditional unit $g$ at the beginning of an OPF period.
$L_n$     Load at bus $n$.
$LimitA_k$     Normal long-term limit of branch $k$ in SCED.
$LimitC_k$     Emergency short-term limit of branch $k$ in SCED.
$ERR_g$     Energy ramp limit for unit $g$ in a SCED interval.
$SRR_g$     Spinning ramp limit for unit $g$ in response to a contingency in 10 minutes.
$x_k$     Reactance of branch $k$.
$OP_g$     Variable operating cost of traditional unit $g$.
$w_s$     Weight (probability) of scenario $s$.
$pIR_{is}^{max}$     Forecasted maximum generation from intermittent renewable unit $i$ in scenario $s$.
$n_{(k-)}$     From-bus of branch $k$.
$n_{(k+)}$     To-bus of branch $k$.
$z^{max}$     Number of lines that can be switched off under an outage.

*Variables*:
$p_{gs}$     Output power of traditional unit $g$ in scenario $s$.
$r_{gs}$     Reserve from traditional unit $g$ in scenario $s$.
$pIR_{is}$     Output power of renewable unit $i$ in scenario $s$.
$cIR_{is}$     Curtailed power of renewable unit $i$ in scenario $s$.
$\theta_{ns}$     Phase angle of bus $n$ in scenario $s$.
$p_{ks}$     Power flow on line $k$ in scenario $s$.
$\theta_{nsc}$     Phase angle of bus $n$ under contingency $c$ in scenario $s$.
$p_{ksc}$     Power flow on line $k$ under contingency $c$ in scenario $s$.
$z_{sck}$     Status of line $k$ under contingency $c$ in scenario $s$.

## I. INTRODUCTION

Renewable generation in the power grid has been increasing significantly in recent years. The fast development of renewables is key to achieving a sustainable power system. High penetration of renewables provides a variety of benefits including energy production cost saving, greenhouse gas emission reduction, and diversification of energy supply resources. However, integration of intermittent renewable power brings a substantial uncertainty to grid operations; moreover, today's industry still uses deterministic optimization to model and solve the optimal power flow (OPF) problem for real-time generation dispatch by assuming intermittent renewable power is perfectly forecasted. Thus, this work examines different stochastic OPF (SOPF) models that can include the probabilistic forecasting information of intermittent renewables and consider multiple scenarios.

Prior efforts in the literature show that network reconfiguration (NR) can provide various benefits to the power system since it is able to reroute power in the transmission network. When NR is modelled as a preventive mechanism, substantial cost savings can be achieved [1]-[4]. However, frequent implementation of network reconfiguration as a preventive control in the base case will significantly accelerate the degradation of circuit breakers and it may also lead to system instability [5]. Thus, NR is often used as a corrective mechanism, which is more realistic and has a potential for industry adoption. It is shown that corrective NR is also able to achieve cost saving [6]-[9], enhance reliability [10]-[14], reduce losses [15]-[18]. Though significant amount of work has been done to demonstrate various benefits provided by NR, its effects on renewable curtailment reduction has not been investigated thoroughly. Renewable generation curtailment is often observed in today's grid operations [19]-[22]. One primary cause of renewable energy curtailment is transmission network congestion. Thus, this work incorporates NR into the SOPF model and shows NR can relieve congestion and reduce congestion-induced undesired renewable curtailment.

In this paper, four different SOPF models are proposed: the relaxed SOPF (R-SOPF), normal SOPF (N-SOPF), enhanced SOPF (E-SOPF), and enhanced SOPF with network reconfiguration (E-SOPFwNR). N-SOPF only includes base-case network constraints while E-SOPF enforces both base-case network constraints and contingency-case network constraints. E-SOPFwNR extends E-SOPF by modelling post-contingency transmission network reconfiguration; R-SOPF that ignores all network constraints is implemented and serves as a benchmark in this work to gauge other SOPF models by comparing their congestion-induced costs. Numerical simulations show that both base-case congestion and contingency-case congestion result in additional system operating cost. By explicitly enforcing contingency-case branch limit constraints in E-SOPF, renewable power curtailment is observed for two out of ten potential scenarios; however, that curtailment can be avoided with the proposed network reconfiguration scheme that is used to relieve transmission congestion in the post-contingency situation. It also shows that NR can substantially reduce congestion cost, especially the contingency-case congestion cost.

The rest of the paper is organized as follows. The proposed methodology and model are presented in Section II. Section III briefly discusses the market implication for the proposed model. Case studies are presented in Section IV. Section V concludes the paper and presents potential future work.

## II. METHODOLOGY AND MODEL

In the real-time operations of electric power systems, system operators must ensure the electric power produced can supply the electric power consumed at all times. To operate a practical power grid, system reliability must be well maintained. Thus, independent system operator (ISO) enforces multiple classes of generating reserves to address system uncertainties such as random load fluctuations and potential outage events. The dispatch solution is required to respect the transmission network capacity such that there are no overloading violations; in addition, *N-1* requirement would force the dispatch solution so that no branch is overloaded even after the system loses one branch. It is also very important to operate the grid economically as a marginal cost saving will be huge due to the large-scale feature of practical power systems. Optimal power flow is the key decision support application for real-time economic dispatch. This section will present the formulation of the stochastic OPF that considers all abovementioned factors and provides solutions for different forecasted scenarios.

*Objective Function*:

The objective of SOPF is presented in (1); it is to minimize the total expected system operating cost.

$$min\ GenCost = \sum_{g \in G} \sum_{s \in S} w_s f(p_{gs}) \quad (1)$$

where $f(p_{gs})$ denotes the cost function for generator $g$. In this work, the linearized cost function in (2) is used for $f(p_{gs})$.

$$f(p_{gs}) = p_{gs} OP_g \quad (2)$$

*Constraints*:

There are a number of different sets of constraints that need to be included in SOPF. Constraint (3) ensures power balance for each bus under each scenario. Generator ramping rate limit is modeled in (4) and output power limit is enforced in (5). In (6), the forecasted maximum generation from intermittent renewables is divided into two parts: scheduled generation and curtailment. Both parts cannot be negative, which is guaranteed by (7). Branch long-term thermal limit constraints are enforced in (8), which ensures there is no transmission violation in the base case. Equation (9) calculates branch flow with two end-bus phase angles. Spinning reserve is modeled in this work and it can only be provided by traditional units that are controllable. Constraints (10) and (11) show that the spinning reserve provided by a unit cannot exceed its ramping capability and available generation capacity respectively. As shown in (12) and (13), the "largest generator" rule is used to set the reserve requirement: the total reserve should be greater than or equal to the largest generation for all scenarios considered in SOPF.

$$\sum_{g \in G(n)} p_{gs} + \sum_{g \in IR(n)} pIR_{gs} + \sum_{k \in K(n+)} p_{ks}$$
$$- \sum_{k \in K(n-)} p_{ks} = L_n \quad \forall n \in N, s \in S \quad (3)$$
$$-ERR_g \leq p_{gs} - p_{g0} \leq ERR_g \quad \forall g \in G, s \in S \quad (4)$$
$$P_{g,min} \leq p_{gs} \leq P_{g,max} \quad \forall g \in G, s \in S \quad (5)$$
$$pIR_{is} + cIR_{is} = pIR_{is}^{max} \quad \forall i \in IR, s \in S \quad (6)$$
$$\{pIR_{is}, cIR_{is}\} \geq 0 \quad \forall i \in IR, s \in S \quad (7)$$
$$-LimitA_k \leq p_{ks} \leq LimitA_k \quad \forall k \in K, s \in S \quad (8)$$
$$p_{ks} = (\theta_{n_{(k-)}s} - \theta_{n_{(k+)}s})/x_k \quad \forall k \in K, s \in S \quad (9)$$
$$0 \leq r_{gs} \leq SRR_g \quad \forall g \in G, s \in S \quad (10)$$
$$p_{gs} + r_{gs} \leq P_{g,max} \quad \forall g \in G, s \in S \quad (11)$$
$$\sum_{m \in G} r_{ms} \geq p_{gs} + r_{gs} \quad \forall g \in G, s \in S \quad (12)$$
$$\sum_{m \in G} r_{ms} \geq pIR_{gs} \quad \forall i \in IR, s \in S \quad (13)$$

Though spinning reserve is enforced in SOPF, reserve deliverability cannot be guaranteed due to post-contingency network congestion. Thus, (14)-(17) should be included in the SOPF model when post-contingency network congestion is a concern. Nodal power balance should be maintained for each contingency case under each system scenario, which is modelled in (14). As shown in (15), for contingency-case network constraints, branch short-term thermal limit is used as the associated flow limit since contingency-case flow is expected to be brought down very soon and it can exceed the normal long-term thermal limit for a short period of time. The relationship between post-contingency line flows and post-contingency phase angles is presented in (16). Equation (17) shows that the flow on the contingency branch is zero.

$$\sum_{g \in G(n)} p_{gs} + \sum_{g \in IR(n)} pIR_{gs} + \sum_{k \in K(n+)} p_{ksc} - \sum_{k \in K(n-)} p_{ksc} = L_n \quad \forall n \in N, c \in C, s \in S \quad (14)$$

$$-LimitC_k \leq p_{ksc} \leq LimitC_k \quad \forall k \in K, c \in C, s \in S \quad (15)$$

$$p_{ksc} = (\theta_{n_{(k-)}sc} - \theta_{n_{(k+)}sc})/x_k \quad \forall k \in \{K - c\}, c \in C, s \in S \quad (16)$$

$$p_{ksc} = 0 \quad \forall k \in \{c\}, c \in C, s \in S \quad (17)$$

As discussed above, network congestion including post-contingency congestion may incur renewable curtailment and network reconfiguration may avoid such undesired curtailment. To include network reconfiguration in the SOPF model, binary variable $z_{sck}$ is introduced to represent the status of switchable line $k$ after contingency $c$ in scenario $s$. Then, (18) can be used to replace (15) and big-M method can be implemented to replace the regular line flow equation (16) with constraints (19)-(20) combined. Constraint (21) shows that the number of lines that are allowed to be switched off is limited to $z^{max}$ that is set to 1 in this paper for stability concern.

$$-z_{sck} LimitC_k \leq p_{ksc} \leq z_{sck} LimitC_k \quad \forall k \in K, c \in C, s \in S \quad (18)$$

$$p_{ksc} - (\theta_{n_{(k-)}sc} - \theta_{n_{(k+)}sc})/x_k + (1 - z_{sck})M \geq 0 \quad \forall k \in \{K - c\}, c \in C, s \in S \quad (19)$$

$$p_{ksc} - (\theta_{n_{(k-)}sc} - \theta_{n_{(k+)}sc})/x_k - (1 - z_{sck})M \leq 0 \quad \forall k \in \{K - c\}, c \in C, s \in S \quad (20)$$

$$\sum_{k \in \{K-c\}} (1 - z_{sck}) \leq z^{max} \quad \forall c \in C, s \in S \quad (21)$$

*SOPF Models*:

The four proposed SOPF models are summarized in Table I. R-SOPF ignores all network constraints and it is implemented to serve as a benchmark to calculate the congestion cost for other SOPF models. As compared to R-SOPF, N-SOPF includes base-case network constraints. E-SOPF extends N-SOPF by enforcing additional contingency-case network constraints. E-SOPFwNR further extends E-SOPF by capturing the flexibility in the transmission network to manage congestion under contingency. Both E-SOPF and E-SOPFwNR can ensure there are no branch overloads and they are *N*-1 secure; however, the solutions obtained with E-SOPFwNR would be more economically efficient.

Table I. Summary of various SOPF models

|  | Power balance constraints | Network constraints | Other constraints |
|---|---|---|---|
| R-SOPF | (3) | N/A | (4)-(7), (9)-(13) |
| N-SOPF | (3) | (8) | (4)-(7), (9)-(13) |
| E-SOPF | (3), (14) | (8), (15) | (4)-(7), (9)-(13), (16)-(17) |
| E-SOPFwNR | (3), (14) | (8), (18) | (4)-(7), (9)-(13), (17), (19)-(21) |

Due to the fact that congestion may exist in the base-case and contingency-case, congestion cost includes two components: base-case congestion cost and contingency-case congestion cost. Since R-SOPF does not include any network constraints and N-SOPF only enforces base-case network constraints, the total congestion cost ($TCC^m$) and total contingency-case congestion cost ($TCCC^m$) for SOPF model $m$ can be defined in

$$TCC^m = TC^m - TC^{R-SOPF} \quad (22)$$
$$TCCC^m = TC^m - TC^{N-SOPF} \quad (23)$$

where $TC^m$, $TC^{R-SOPF}$, and $TC^{N-SOPF}$ denote the total operating cost for SOPF model $m$, R-SOPF, and N-SOPF respectively.

### III. MARKET IMPLICATION

Two-thirds of the U.S.'s electricity demand is served in the territory of wholesale electric energy markets. Locational marginal pricing is the essential pricing mechanism for clearing the wholesale energy markets. The electricity prices, often referred to as locational marginal prices (LMP), at different locations could be very different due to network congestion and losses. In real-time markets, LMP is determined along with generator dispatch points by solving OPF. Thus, it is very important to analyze the impact of including network reconfiguration in the model on the market.

For R-SOPF and N-SOPF, there are only base-case power balance constraints that involve bus load. Thus, the nodal LMP for R-SOPF and N-SOPF can be determined as follows,

$$LMP_n = \sum_{s \in S} \delta_{ns} \quad (24)$$

where $\delta_{ns}$ denotes the dual of constraint (3) for bus $n$ in scenario $s$.

For E-SOPF and E-SOPFwNR, in addition to base-case power balance constraints, they also include post-contingency power balance constraints that involve bus load. Thus, the nodal LMP for E-SOPF and E-SOPFwNR can be determined as follows,

$$LMP_n = \sum_{s \in S} \delta_{ns} + \sum_{s \in S} \sum_{c \in C} \delta_{nsc} \quad (25)$$

where $\delta_{nsc}$ denotes the dual of constraint (14) for bus $n$ in scenario $s$ under contingency $c$.

In addition to nodal LMP, the system-wide average LMP ($AvgLMP$) and weighted LMP ($AvgLMP^w$) are defined as follows,

$$AvgLMP = \sum_{n \in N} LMP_n \quad (26)$$
$$AvgLMP^w = \sum_{n \in N} LMP_n L_n / \sum_{n \in N} L_n \quad (27)$$

Equations (28), (29) and (30) define the load payment, traditional generator revenue and intermittent renewables revenue respectively. Traditional generator profit is calculated by (31). Congestion revenue, which is different with congestion cost, can be calculated by (32).

$$LdPaymt = \sum_{n \in N} (LMP_n L_n) \quad (28)$$
$$GenRvn = \sum_{n \in N} (LMP_n (\sum_{g \in G(n)} \sum_{s \in S} w_s p_{gs})) \quad (29)$$
$$ResGenRvn = \sum_{n \in N} (LMP_n (\sum_{g \in IR(n)} \sum_{s \in S} w_s \cdot pIR_{gs})) \quad (30)$$
$$GenProfit = GenRvn - GenCost \quad (31)$$
$$CongRvn = LdPaymt - GenRvn - ResGenRvn \quad (32)$$

### IV. CASE STUDIES

The proposed four SOPF models are tested on a modified area of the IEEE RTS-96 test system. This one area test case has 24 buses, 38 branches, 33 traditional generators and 5 intermittent generators (IG). 22 traditional generators are synchronized to the grid, while the 11 offline traditional generators are considered as unavailable in SOPF in this work. Loads are connected to 17 buses and correspond to a total demand of 2,850 MW. The load profile is shown in Fig. 1. The network topology for this test system is the same with the case shown in [23]. The 5 intermittent generators are located at bus 1, 14, 15, 21, and 22 respectively and their total production is 650 MW. Ten different scenarios of intermittent renewable generation are considered in this work and they are shown in Fig. 2.

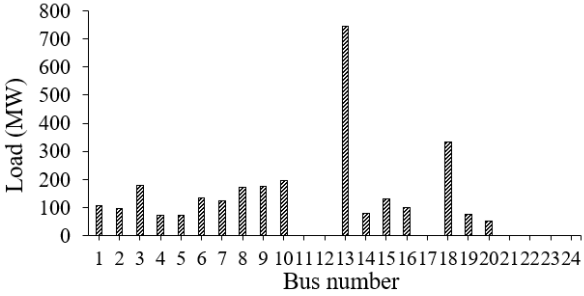
Fig. 1. Load profile of the IEEE 24-bus test system.

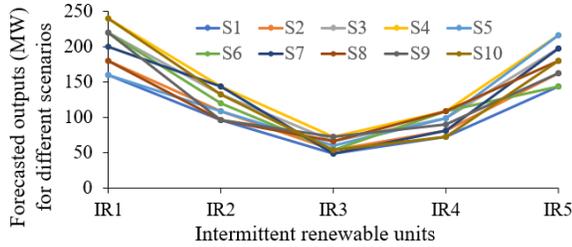
Fig. 2. Renewable generation forecasted outputs under different scenarios.

Table II and Fig. 3 present the cost information for the four proposed SOPF models. As expected, the total cost of R-SOPF is minimum as its network is assumed to have infinite capacity while E-SOPF corresponds to the highest cost as it includes both base-case network constraints and contingency-case network constraints. It is also observed that including network reconfiguration in E-SOPFwNR can achieve significant cost saving, especially congestion cost reduction. Table II shows that including NR can achieve a reduction of 34.8% in total contingency-case congestion cost. This indicates that NR can effectively relieve network congestion in post-contingency situations. The solutions obtained with both E-SOPF and E-SOPFwNR will not lead to any violation concerns; however, the solution obtained with E-SOPFwNR is more attractive as its cost is lower. Moreover, with E-SOPFwNR, there is no renewable curtailment under any scenarios. However, with E-SOPF, renewable generation curtailment is forced in two scenarios (4 & 5) to maintain system reliability, which corresponds to a power curtailment of 14.2 MW and 3.8 MW respectively.

Including additional variables and constraints will increase the model computational complexity and it is expected that the proposed E-SOPFwNR model is more complex to solve as compared to the other three models, R-SOPF, N-SOPF and E-SOPF. This is demonstrated by the results presented in Table III. Though E-SOPFwNR can lower the cost and reduce renewable energy curtailment, it takes much longer time to solve. Therefore, advanced decomposition algorithm is required to accelerate the solution process, which would be the future work of this paper.

Table II. Various costs for different SOPF models

|  | R-SOPF | N-SOPF | E-SOPF | E-SOPFwNR |
|---|---|---|---|---|
| $TC^m$ ($/h) | 36,346 | 39,536 | 44,965 | 43,075 |
| $TCC^m$ ($/h) | 0 | 3,190 | 8,620 | 6,729 |
| $TCCC^m$ ($/h) | 0 | 0 | 5,430 | 3,540 |
| Reduction in congestion cost with NR ($/h) | N/A | N/A | N/A | 1,890 (21.9%) |
| Reduction in contingency-case congestion cost with NR ($/h) | N/A | N/A | N/A | 1,890 (34.8%) |

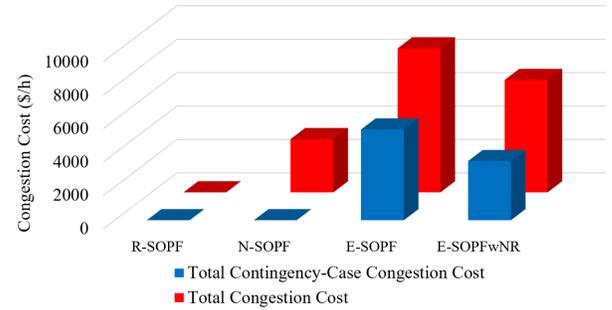
Fig. 3. Congestion cost for various SOPF models.

Table III. Model complexity statistics after AMPL Presolve Process

|  | R-SOPF | N-SOPF | E-SOPF | E-SOPFwNR |
|---|---|---|---|---|
| # of variables | 1,110 | 1,110 | 23,600 | 37,190 |
| # of constraints | 1,040 | 1,040 | 23,530 | 64,730 |
| # of nonzeros | 8,460 | 8,460 | 86,660 | 222,700 |
| Solution time (s) | 0.016 | 0.032 | 0.344 | 144.938 |

Table IV. Market results with various SOPF models

|  | R-SOPF | N-SOPF | E-SOPF | E-SOPFwNR |
|---|---|---|---|---|
| $AvgLMP$ ($/MWh) | 39.9 | 39.8 | 46.5 | 32.5 |
| $AvgLMP^w$ ($/MWh) | 39.9 | 39.6 | 48.2 | 35.5 |
| $LdPaymt$ ($/h) | 113,715 | 112,854 | 137,385 | 101,299 |
| $ResGenRvn$ ($/h) | 25,967 | 18,426 | 13,209 | 11,446 |
| $GenRvn$ ($/h) | 87,748 | 62,989 | 72,706 | 61,884 |
| $GenProfit$ ($/h) | 42,441 | 14,493 | 18,780 | 9,848 |
| $CongRvn$ ($/h) | 0.0 | 31439 | 51,471 | 27,969 |

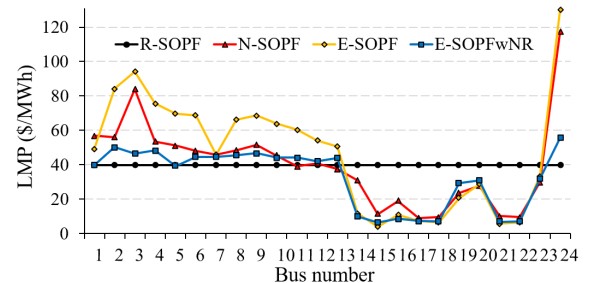
Fig. 4. Nodal LMP for various SOPF models.

Table IV presents the market results with various SOPF models. When network reconfiguration is considered in SOPF, the load payment drops significantly, as well as the traditional generator revenue and profit, renewable generator revenue and congestion revenue. Fig. 4 illustrates nodal LMP with various SOPF models proposed in this work. The LMPs at different buses for R-SOPF are the same throughout the entire network since losses are not explicitly modelled in this work and R-SOPF assumes infinite capacity network. Though the nodal LMP profiles obtained with N-SOPF, E-SOPF, and E-SOPFwNR are variable due to congestion, they share similar trends. Overall, the LMPs for E-SOPFwNR are lower than E-SOPF, which indicates that NR can improve social welfare.

## V. Conclusion and Future Work

Fast development of renewable energy is key to achieving the next generation smart grid. However, due to network congestion, undesired renewable power curtailment is often observed, which is a waste of resources and discourages further deployment of renewables. Today's operational tool still uses deterministic OPF that models the transmission network as static asset. Thus, to capture the probabilistic forecasting information of intermittent renewables and utilize network reconfiguration in grid real-time operations, we propose an enhanced SOPF model with consideration of network reconfiguration in this paper. Numerical simulations show that post-contingency network constraints may lead to renewable generation curtailment; moreover, including network reconfiguration in SOPF can avoid such undesired curtailment and in the meantime, it significantly reduces congestion cost.

Though the proposed E-SOPFwNR can utilize the network reconfiguration scheme to relieve network congestion and facilitate grid integration of renewables, it is much more computationally complex than the current deterministic OPF model. This means it will take a lot more time to solve E-SOPFwNR for large-scale practical systems, which makes it hard for industry adoption. Thus, future work will involve developing effective algorithms for solving large-scale stochastic optimal power flow problems that incorporate contingency-case network constraints and network reconfiguration scheme.